\documentclass[12pt,leqno]{article}
\usepackage{amsmath,amsfonts,amscd,ProtocoloS3H,url,txfonts}
\usepackage{pdfsync}
\usepackage{graphicx}
\usepackage{color}
\usepackage[all]{xy}

\title{\Large\bf
{The Gysin sequence \\ for $\S^{3}$-actions on
manifolds}\footnote{This work has been partially supported 
by the UPV/EHU grant EHU09/04 and by the Spanish MICINN grant MTM2010-15471.} 
\color{white}.\normalcolor\footnote{Key Words and Phrases: Gysin sequence, basic
cohomology, $\S^3$-actions, controlled forms.} 
\color{white}.\normalcolor\footnote{2010 Mathematical Subject
Classification. Primary 57R19; Secondary 57R30, 57S15.}
}
\date{}

\author{
Jos\'e Ignacio Royo Prieto\thanks{Departamento de Matem\'atica
Aplicada, 
University of the Basque Country UPV/EHU. Alameda de Urquijo
s/n. 48013 Bilbao, SPAIN. \hspace{5cm}  {\sl joseignacio.royo@ehu.es}.}
\\ 
\\
{\small University of the Basque Country UPV/EHU}
\and Martintxo
Saralegi-Aranguren\thanks{Universit\'e d'Artois, Laboratoire de
Math\'ematiques de Lens EA~2462. F\'ed\'eration CNRS
Nord-Pas-de-Calais FR~2956. Facult\'e des Sciences Jean Perrin. Rue
Jean Souvraz, S.P.\kern 1mm 18.
 F-62\kern 1mm 300 Lens, France.
{\sl saralegi@euler.univ-artois.fr}. }
\\
\\ {\small Universit\'e d'Artois }
}

\begin{document}

\maketitle

\thispagestyle{empty}
\begin{abstract}
Given a smooth action of $\bi{\S}{3}$ on a  manifold $M$, we are
interested in the relationship between the cohomologies of $M$ and
$M/\bi{\S}{3}$. If the action is free, we have indeed a principal
$\bi{\S}{3}$-bundle, and this relationship is described by the
classical Gysin sequence, which also exists when the action is
semi-free (i.e., fixed points are allowed) \cite{Br}. In this work,
we obtain a Gysin sequence for the case of a general smooth action.
An exotic term appears, and we show that it is an obstruction for
the duality of the second term of the de Rham spectral sequence
associated to the action.
\end{abstract}

\bigskip

Let us consider a smooth action $\Phi \colon G \times M \to
M$ of a compact Lie group on a manifold $M$. The action $\Phi$ induces naturally a filtration
$\{ \bi{F}{i}\lau{\Om}{*}{}{M}\tq i \in \N \}$ of the complex of de Rham differential forms $\lau{\Om}{*}{}{M}$, defined by:
\begin{equation*}
\label{filtra}
\bi{F}{i}\lau{\Om}{i+j}{}{M}  = \left\{ \om \in \lau{\Om}{j}{}{M} \tq
i_{X_0} \cdots i_{X_{j}}\om= 0
\hbox{ for each family } \left\{ X_0, \ldots ,
X_{j}\right\} \subset \homo{\mathfrak{X}}{\Phi}{M}
\right \}.
\end{equation*}
Here, we have denoted by $\homo{\mathfrak{X}}{\Phi}{M}$
the orbit distribution of $TM$ formed by the vector fields of $M$ tangent to the orbits of $\Phi$.
This filtration defines the first quadrant {\em de Rham spectral sequence}, which converges to $\lau{H}{*}{}{M}$. The underlying motivation of this paper is the study of the Poincar\'e duality  of the second term $\hiru{E}{s,t}{2}$ of this spectral sequence.

When $\Phi$ is free, we have the duality $\hiru{E}{s,t}{2} \cong
\hiru{E}{n-s,\ell-t}{2}$, where $n = \dim M/G$ and $\ell = \dim G$. This property
is lost when the action is no longer free.

Inspired by the work of Goresky and MacPherson, one expects to recover the Poincar\'e duality by using intersection cohomology. This is the case when the group $G$ is the circle $\sbat$  (see \cite{JI}).
The next natural group $G$ to study is $\bi{\S}{3}$ (of rank 1 and not abelian). It has been proved in \cite{S3} that Poincar\'e duality  still holds when the action $\Phi$ is semi-free.

What about the other $\bi{\S}{3}$-actions? Surprisingly, the second term of the above spectral sequence has not been computed yet in this context. The main result of this paper is
in the following Gysin sequence, which computes this second term:

\begin{equation*}
\xymatrix @R=4pc{
\cdots \ar[r] &
\coho{H}{i}{M}\ar[r]^-{\footnotesize \textcircled{\tiny  2}}
&\smash{\underbrace{
 \coho{H}{i-3}{  M/\bi{\S}{3},\Sigma/\bi{\S}{3}  }
  }_{  \hiru{E}{i-3,3}{2}  }}
 \oplus
\smash{\underbrace{ \left( \coho{H}{i-2}{M^{\sbat}} \right)^{-\ib{\Z}{2}}}_{\hiru{E}{i-2,2}{2}}}
\ar[r]^-{\footnotesize \textcircled{\tiny  3}}
&
\smash{\underbrace{\coho{H}{i+1}{M/\bi{\S}{3}}}_{\hiru{E}{i+1,0}{2}} \ar[r]^-{\footnotesize \textcircled{\tiny  1}}}
 &  \coho{H}{i+1}{M}\ar[r] & \cdots &
}
\end{equation*}
with $\hiru{E}{i,1}{2} =0$, where
\begin{itemize}
\item[-]$\Sigma$ is the subset of  points of $M$ whose  isotropy group is infinite;

\item[-] our choice of maximal torus, $\S^1$, is $\left\{a+bi\tq a^2+b^2=1\right\}\le\S^3$;

\item[-] the $\ib{\Z}{2}$-action is induced by $j \in \bi{\S}{3}$,

\item[-] $(-)^{-\ib{\Z}{2}}$ denotes the subspace of antisymmetric elements (cf. \refp{anti}),

\item[-] {\footnotesize \textcircled{\tiny  1}} is induced
by the natural projection $\pi \colon M \to M/\bi{\S}{3}$,

\item[-] {\footnotesize \textcircled{\tiny  2}} is induced by the
integration along the fibers of $\pi $, and

\item[-]{\footnotesize \textcircled{\tiny  3}} involves the multiplication by the
{\em Euler class}  $[e] \in \lau{\IH}{4}{\per{4}}{M/\bi{\S}{3}}$
(cf. \cite{S3})

\end{itemize}
 (cf. Theorem \ref{gysin} and paragraph 2.5).

Notice that the first floor  $\hiru{E}{i,1}{2}$
always vanishes (even for any perversity!)  whereas the second floor
$\hiru{E}{i,2}{2}$ may not, as we show in example 2.4. So, it follows   that Poincar\'e duality does not work
in the generic case, even considering intersection cohomology.

As a consequence, a new approach is needed in order to extend the Poincar\'e duality  for general actions, and, thus,
for Singular Riemannian Foliations, as is the partition induced by the orbits of a general $\mathbb{S}^3$-action. Results in this direction are being explored.

In fact, in \cite{Mas} the duality of the spectral sequence associated to a Singular Riemannian Foliation
is claimed for the special case of extreme perversities (which is tantamount to working only in the
regular stratum or relatively to the singular strata). As the previous counterexample shows, a new approach is needed
in order to extend the  duality result for general perversities if one wants to consider the usual cohomology ($\per{p}=\per{0}$) case.

A different Gysin sequence relating the cohomology of $M$ and the $\S^3$-equivariant cohomology of $M$ was constructed in \cite{chen}, also in the case of a general smooth $\S^3$-action.

\medskip

In the sequel $M$  is a connected, second countable, Hausdorff,
without boundary and smooth
(of class $C^\infty$) manifold.
We fix a smooth action $\Phi \colon \bi{\S}{3} \times M \to M$.

\medskip

We wish to thank the referees for the indications given in order to improve this paper.

\prg{Stratifications and differential forms}

We describe the stratification arising from the action. We also
introduce the controlled differential forms, defined by Verona, in
order to compute the singular cohomology in this context.

\prgg{Thom-Mather structure}

There are three possibilities for the dimension of the isotropy
subgroup\footnote{We refer the reader to \cite{Br} for the notions
related with compact Lie group actions, such as isotropy, invariant
tubular neighborhoods,\ldots }
 $\bi{\S}{3}_x$ of a point $x \in M$, namely: 0,1 and 3.
 So, we have the dimension-type filtration
\begin{equation*}
 F = \left\{ x \in M \tq  \dim \bi{\S}{3}_x  = 3\right\} \subset \Sigma = \left\{ x \in M \tq  \dim \bi{\S}{3}_x  \geq 1\right\}\subset M=\left\{ x \in M \tq  \dim \bi{\S}{3}_x  \geq 0\right\} .
\end{equation*}
In this section, we describe the geometry of the triple
$(M,\Sigma,F)$. The subset $\Sigma$ is not necessarily a manifold,
but  the subsets $F = M^{\bi{\S}{3}}$, $\Sigma \menos F=\left\{ x \in M
\tq  \dim \bi{\S}{3}_x  = 1\right\}$ and $M\menos \Sigma = \left\{ x
\in M \tq  \dim \bi{\S}{3}_x  = 0\right\}$ are  proper invariant
submanifolds\footnote{In fact, these manifolds may have  connected
components with different dimensions.} of $M$. So, we can consider
$\tau_0 \colon T_0\to F$ and  $\tau_1 \colon T_1\to \Sigma\menos F$
two invariant  tubular neighborhoods in $M$. Over each connected
component, the  structure group is the orthogonal group.  Associated
to these tubular neighborhoods we have the following maps ($k=0,1$):

\medskip

$\rightsquigarrow$ The {\em radius map}
$\nu_k \colon T_k \to [0,\infty[$, defined fiberwise by $u
\mapsto \| u\|$. It is an invariant  smooth map.

\medskip

$\rightsquigarrow$
The {\em dilatation map} $\partial_k \colon [0,\infty[ \times T_k \to T_k$,
defined fiberwise by $(t,u) \mapsto t \cdot u$.  It is a smooth  equivariant map.

The family of tubular neighborhoods $\mathfrak{T} _{M}= \{T_0,T_1\}$ is
a {\em Thom-Mather system} when:

\medskip

(TM)
   $\left\{
   \begin{array}{c}    \tau_{0} = \tau_{0} \rondp \tau_{1}\\   \nu_{0} = \nu_{0} \rondp \tau_{1}  \end{array} \right\} $ on
$T_{0} \cap T_{1} = \tau_{1}^{-1}(T_{0} \cap
   (\Sigma\menos F))$.

\bl \label{ThM} Thom-Mather systems exist. \el \pro We fix an
invariant tubular neighborhood ${\tau}_{0} \colon T_0 \to F$. It
exists since $F$ is an invariant closed submanifold of $M$. Since
the isotropy subgroup of any point of $F$ is the whole $\bi{\S}{3}$,
we can find\footnote{For each connected component of $F$.} an atlas
$\mathcal{\mathcal{A}} = \left\{ \phii \colon U \times \bi{\R}{n}
\to \tau_0^{-1}(U)\right\}$ of $\tau_0$, having $O(n)$ as structure
group,  and an orthogonal action $\Psi \colon \bi{\S}{3} \times \R^n
\to \R^n$ such that \be \label{gs} \phii (x, \Psi(g,v)) =
\Phi(g,\phii(x,v) )\quad \forall x \in U, \forall v \in \R^n \hbox{
and } \forall g \in \bi{\S}{3}. \ee We write ${\tau}_{0}' \colon S_0
\to F$ the restriction of $\tau_0$, where $S_0$ is the submanifold
$\nu^{-1}_0(1)$. It is a fiber bundle. The restriction ${\tau}_{0}''
\colon (S_0  \cap (\Sigma \menos F)) \to F$ is also a fiber bundle
whose induced atlas is $\mathcal{\mathcal{A}}'' = \left\{ \phii
\colon U \times \S_{\Sigma}^{^{n-1}}\to {\tau''_0}^{-1}(U)\right\}$,
where $\S_{\Sigma}^{^{n-1}} = \left\{ w \in \bi{\S}{n-1} \tq  \dim
\bi{\S}{3}_w  = 1\right\}$.

The map $\mathfrak{L}_0 \colon T_0\menos F \to S_0 \times
]0,\infty[$, defined by $\mathfrak{L}_0(x) = \left( \partial_0\left(
\nu_0(x)^{-1},x\right),\nu_0(x)\right)$, is an equivariant
diffeomorphism. Under $\mathfrak{L}_0$:

$\rightsquigarrow$ the map $\tau_0$ becomes $(y,t) \mapsto \tau_0'(y)$,

\smallskip

$\rightsquigarrow$
 the map $\nu_0$ becomes $(y,t) \mapsto t$, and

\smallskip

$\rightsquigarrow$
the manifold $T_0 \cap (\Sigma \menos F)$ becomes $(S_0 \cap (\Sigma \menos F)) \times ]0,\infty[$.

Since the structure group  of  ${\tau}_{0}' $ is a  compact Lie group, condition \refp{gs} allows us to construct an invariant Riemannian metric $\mu_0$ on $S_0$ such that the fibers of $\tau_0'$ are totally geodesic submanifolds and $\left( T(S_0 \cap (\Sigma \menos F)) \right)^{\bot} \subset \ker \left( \tau'_0\right)_*$. Then, if we consider the associated tubular neighborhood $\tau_1' \colon T_1'\to   S_0 \cap (\Sigma \menos F)$ we have $\tau_0' \circ \tau_1' = \tau_0'$.

We can construct now an invariant Riemannian metric $\mu$ on $M\menos F$ such that under $\mathfrak{L}_0$:

$\rightsquigarrow$ the metric $\mu$ becomes $\mu_0 +dr^2$ on $S_0 \times ]0,\infty[$.

\nt  We consider the associated tubular neighborhood $\tau_1 \colon T_1 \to \Sigma\menos F$. Verification of the property (TM) must be done on $T_0 \cap T_1$, where using $\mathfrak{L}_0$, we get:

$\rightsquigarrow$ $T_0 \cap T_1$ becomes $T_1' \times ]0,\infty[$.

\smallskip

$\rightsquigarrow$ $\tau_1$ becomes $(y,t) \mapsto (\tau'_1(y),t)$.

\nt A straightforward calculation gives (TM) and ends the proof. \qed

We fix a such system $\mathfrak{T}_M$. For each $k \in \{0,1\}$, we shall write $D_{k} \subset M$ the open
subset  $ \nu_{k}^{-1}([0,1[)$ and call it the {\em soul} of the tubular neighborhood $\tau_k$. We shall write $\Delta_0 = D_0 \cap \Sigma$.

\prgg{Verona's differential forms} As it is shown in \cite{V3},
 the singular cohomology of $M$ can be computed
by using differential forms on $M\menos \Sigma$.
This is the tool we use in this work.
The complex of {\em controlled forms} (or {\em Verona's forms})
 of $M$  is defined by
\bee
\lau{\Om}{*}{V}{M}=\left\{ \om \in \coho{\Om}{*}{M \menos \Sigma } \tq \exists
\om_{1} \in \lau{\Om}{*}{}{\Sigma \menos F}
\hbox{ and }
\om_{0} \in \coho{\Om}{*}{F} \hbox{ with }
\left\{
\begin{array}{l}
(a) \ \tau_{1}^* \om_{1} = \omega \ \hbox{ on } D_{1}  \menos\Sigma \\
(b) \ \tau_{0}^* \om_{0} = \omega \ \hbox{ on } D_{0}  \menos \Sigma  \\
(c) \ \tau_{0}^* \om_{0} = \omega_1 \ \hbox{ on } \Delta_{0}  \menos F
\end{array}
\right\}
\right\}.
\eee
Following \cite{V3} we know that the cohomology  of the complex $\lau{\Om}{*}{V}{M}$ is the singular cohomology $\coho{H}{*}{M}$.

We also use in this work the complex $\lau{\Om}{*}{V}{\Sigma}= \left\{ \gamma \in \coho{\Om}{*}{\Sigma \menos F} \tq \exists \gamma_{0} \in \lau{\Om}{*}{}{F} \hbox{ with }
\tau_{0}^* \gamma_{0} = \gamma\ \hbox{ on } \Delta_0\menos F\right\}$ and
the {\em relative complexes}
$\lau{\Om}{*}{V}{M,\Sigma}= \{ \om \in \lau{\Om}{*}{V}{M} \tq \om_{1} \equiv 0\}$ and
$\lau{\Om}{*}{V}{\Sigma,F}= \{ \gamma \in \lau{\Om}{*}{V}{\Sigma } \tq  \gamma_{0} \equiv 0\}$.

\smallskip

Since $M$ is a manifold, controlled forms are in fact differential forms on $M$.
\bl
\label{ext}
Any controlled form of $M$ is the restriction of a differential form of $M$.
\el
\pro
First, we construct  a section $\sigma$ of the restriction $\rho \colon \lau{\Om}{*}{V}{M} \to \lau{\Om}{*}{V}{\Sigma}$ defined by $\rho (\omega) = \omega_1$. Let us consider a smooth function $f \colon ]0,\infty[ \to [0,1]$ verifying $f \equiv 0$ on $[3,\infty[$ and $f \equiv 1$ on $]0,2]$. Notice that the compositions $ f  \rondp  \nu_0 \colon M\to [0,1]$ and $ f  \rondp  \nu_1 \colon M\menos F\to [0,1]$ are smooth invariant maps.
So, for each $\gamma \in \lau{\Om}{*}{V}{\Sigma}$ we have
\begin{equation}
\label{sigma}
\sigma (\gamma) = ( f  \rondp  \nu_0)  \cdot \tau^*_0 \gamma_0 + (1 - ( f  \rondp  \nu_0)) \cdot  ( f  \rondp  \nu_1) \tau^*_1 \gamma\in \lau{\Om}{*}{}{M}.
\end{equation}
This differential form is a controlled form since
\begin{itemize}
\item[]
\begin{itemize}
\item[(a)] Since $(f \rondp \nu_1) \equiv 1 $ on $D_1$,  $(f \rondp \nu_0) \equiv 0 $ on $M \menos T_0$ and (TM) then we have
\bee
\sigma (\gamma)  =
 ( f  \rondp  \nu_0) \cdot \tau^*_1\tau^*_0 \gamma_0+ (1 - ( f  \rondp  \nu_0)) \cdot   \tau^*_1 \gamma
=
 \tau^*_1  \left( ( f  \rondp  \nu_0) \cdot\tau^*_0 \gamma_0+ (1 - ( f  \rondp  \nu_0)) \cdot  \gamma \right)
  \eee
 on $ D_{1} \menos \Sigma$. This gives $(\sigma (\gamma))_1 =  ( f  \rondp  \nu_0) \cdot\tau^*_0 \gamma_0+ (1 - ( f  \rondp  \nu_0)) \cdot  \gamma$. Since $\tau^*_0 \gamma_0 = \gamma$ on $\Delta_0 \menos F$ then $(\sigma (\gamma))_1 =  ( f  \rondp  \nu_0) \cdot \gamma+ (1 - ( f  \rondp  \nu_0)) \cdot  \gamma = \gamma$.
\item[(b)]  Since $(f \rondp \nu_0) \equiv 1 $ on $D_0$ then we have
$ \sigma (\gamma) = \tau^*_0 \gamma_0   $ on $ D_{0} \menos \Sigma$. This gives  $(\sigma (\gamma))_0 = \gamma_0$.
\item[(c)] We have
$(\sigma (\gamma))_1 =  \gamma =\tau^*_0  \gamma_0 =\tau^*_0 (\sigma (\gamma))_0 $  on $ \Delta_{0} \menos F$.
\end{itemize}
\end{itemize}
This map  $\sigma$ is a section of $\rho$ since $\rho (\sigma (\gamma)) = (\sigma (\gamma))_1 = \gamma.$

In particular, $\rho (\om -  \sigma(\rho (\om )))=0$ for each $\om \in \lau{\Om}{*}{V}{M}$.
As $\sigma(\rho(\omega))\in\lau{\Om}{*}{}{M}$ (cf. \refp{sigma}) and coincides with $\omega$ in the open set $(D_0 \cup D_1)\menos \Sigma$ we conclude that $\omega$  can be extended to $M$.
\qed

\prgg{Invariant forms}

We fix $\{u_1,u_2,u_3\}$
a basis of the Lie algebra of  $\bi{\S}{3}$  with $[u_1,u_2] = u_3$, $[u_2,u_3] = u_{1}$ and
$[u_3,u_1] = u_{2}$.
We denote by $X_{i} \in \homo{\mathfrak{X}}{\Phi}{ M}$ the {\em fundamental vector
field} associated
to $u_i$, $i=1,2,3$.

A controlled form $\om$ of $M$ is an {\em invariant form} when $\ib{L}{X_i}\om  =0$
for each  $i=1,2,3$.
The complex of invariant forms is denoted by
$\lau{\BOm}{*}{V}{M} $. The inclusion $\lau{\BOm}{*}{V}{M} \hookrightarrow \lau{\Om}{*}{V}{M} $ induces an isomorphism in cohomology. This a standard argument based on the fact that $\bi{\S}{3}$ is a connected compact Lie group (cf. \cite[Theorem I, Ch. IV, vol. II]{GHV}). So,
\begin{equation}
\label{ver}
\coho{H}{*}{\lau{\BOm}{.}{V}{M}} = \coho{H}{*}{\lau{\Om}{.}{V}{M}} =\coho{H}{*}{M}.
\end{equation}

\prgg{Basic forms} A controlled form $\om$ of $M$ is a {\em basic form} when $\ib{i}{X}\om =  \ib{i}{X}d\om =0$ for each  $X \in \homo{\mathfrak{X}}{\Phi}{M}$.
The complex of the basic forms is denoted by $
\lau{\Om}{*}{V}{M/\bi{\S}{3}}$. In a similar fashion we define $
\lau{\Om}{*}{V}{\Sigma/\bi{\S}{3}}$.
In this work, we shall use the following relative versions of these complexes:
$\lau{\Om}{*}{V}{M/\bi{\S}{3},\Sigma/\bi{\S}{3}} = \lau{\Om}{*}{V}{M/\bi{\S}{3}} \cap \lau{\Om}{*}{V}{M,\Sigma}$, as well as $\lau{\Om}{*}{V}{\Sigma/\bi{\S}{3}, F}=\lau{\Om}{*}{V}{\Sigma/\bi{\S}{3}}\cap\lau{\Om}{*}{V}{\Sigma,F}$.
\bl
\label{rela}
\begin{equation*}
\lau{H}{*}{}{\lau{\Om}{\cdot}{V}{M/\bi{\S}{3}}} = \coho{H}{*}{M/\bi{\S}{3}} \hspace{1cm} \hbox{ and } \hspace{1cm}
\lau{H}{*}{}{\lau{\Om}{\cdot}{V}{M/\bi{\S}{3},\Sigma/\bi{\S}{3}}} =
\coho{H}{*}{M/\bi{\S}{3},\Sigma/\bi{\S}{3}}.
\end{equation*}
\el
\pro The orbit space $M/\bi{\S}{3}$ is a stratified pseudomanifold.
The family of tubular neighborhoods
$\mathfrak{T}_{M/\bi{\S}{3}}= \{\pi (T_0),\pi ( T_1)\}$ is
a {\em Thom-Mather system}. Here, $\pi \colon M \to  M/\bi{\S}{3}$ denotes
the canonical projection. Using this projection, we identify the complex of controlled forms
 of $  M/\bi{\S}{3}$ with $\lau{\Om}{\cdot}{V}{M/\bi{\S}{3}}$, and the same holds for
$\Sigma$.

Since
$\lau{H}{*}{}{\lau{\Om}{\cdot}{V}{X}} = \coho{H}{*}{X}$ for any stratified pseudomanifold $X$,
then  $\lau{H}{*}{}{\lau{\Om}{\cdot}{V}{M/\bi{\S}{3}}} = \coho{H}{*}{M/\bi{\S}{3}}$ and  $\lau{H}{*}{}{\lau{\Om}{\cdot}{V}{\Sigma/\bi{\S}{3}}} = \coho{H}{*}{\Sigma/\bi{\S}{3}}$
 (cf. \cite{V3}).
In fact, the orbit spaces ${M/\bi{\S}{3}}$ and $\Sigma / \bi{\S}{3}$ are triangulable \cite{V}, and by \cite{Y}, both of them possess good coverings. Moreover, any open covering of  ${M/\bi{\S}{3}}$ (resp. $\Sigma / \bi{\S}{3}$) possesses a subordinated partition of unity made up of controlled functions.
 So, we can proceed as in \cite{BT} and construct a
commutative diagram
\bee
\xymatrix{
\cdots \ar[r]   &
\coho{H}{p}{\lau{\Om}{\cdot}{V}{M/\bi{\S}{3},\Sigma/\bi{\S}{3}}}\ar[r]&
\coho{H}{p}{\lau{\Om}{\cdot}{V}{M/\bi{\S}{3}}}\ar[r] \ar[d]^{\ib{f}{M}}&
\coho{H}{p}{\lau{\Om}{\cdot}{V}{\Sigma/\bi{\S}{3}}}\ar[r]\ar[d]^{\ib{f}{\Sigma}}&
\coho{H}{p+1}{\lau{\Om}{\cdot}{V}{M/\bi{\S}{3},\Sigma/\bi{\S}{3}}}\ar[r]&\cdots\\
\cdots \ar[r]   &
\coho{H}{p}{M/\bi{\S}{3},\Sigma/\bi{\S}{3}} \ar[r]&
\coho{H}{p}{M/\bi{\S}{3}}\ar[r]  &
\coho{H}{p}{\Sigma/\bi{\S}{3}} \ar[r]\ar[r]&
  \coho{H}{p+1}{M/\bi{\S}{3},\Sigma/\bi{\S}{3}}\ar[r]&
\cdots
}
\eee
where the vertical arrows are isomorphisms and the horizontal rows are the long  exact sequences associated to the pair $(M/\bi{\S}{3},\Sigma/\bi{\S}{3})$. This gives $\lau{H}{*}{}{\lau{\Om}{\cdot}{V}{M/\bi{\S}{3},\Sigma/\bi{\S}{3}}} = \coho{H}{*}{M/\bi{\S}{3},\Sigma/\bi{\S}{3}}$\footnote{Notice that this is not the five lemma.}.
\qed

\prg{Gysin sequence}

We construct the long exact sequence associated to the action $\Phi \colon \bi{\S}{3} \times M \to M$ relating the cohomology of $M$ and $M/\bi{\S}{3}$.
First of all, we shall use strongly that $\Phi$ is {\em almost free}\footnote{All the isotropy subgroups are
finite groups.} in $ M\menos \Sigma$ to get a better description of the controlled forms of $M$.

\prgg{Decomposition of a differential form}

We endow $ M\menos \Sigma$  with an $\bi{\S}{3}$-invariant Riemannian metric
$\mu_0$, which
exists because $\S^{3}$ is compact.  We also fix a
bi-invariant Riemannian metric $\nu$ on the Lie group $\bi{\S}{3}$.
Consider  now the $\mu_0$-orthogonal
$\bi{\S}{3}$-invariant decomposition
$
T (M\menos \Sigma) =  \mathcal{D} \oplus \xi,
$
where $\mathcal{D}$ is the distribution generated by $\Phi$.
Since the action  $\Phi$  is almost free  on $M\menos \Sigma$, for each point $x \in  M\menos \Sigma$, the  family $\{ X_1(x), X_2(x), X_3(x)\}$ is a basis of $\mathcal{D}_x$.
We define the $\bi{\S}{3}$-Riemannian metric $\mu$ on $ M\menos \Sigma$ by putting
$$
\mu (w_1, w_2) =
\left\{
\begin{array}{ll}
\mu_0(w_1,w_2) & \hbox{if } w_1,w_2 \in \xi_x \\
0 & \hbox{if } w_1 \in \xi_x, w_2 \in\mathcal{D}_x \\
\delta_{i,j}& \hbox{if } w_1 = X_i(x), w_2 = X_j(x)
\end{array}
\right.
$$

We denote by $\chii_{i} = i_{X_i} \mu\in \lau{\Om}{1}{}{ M\menos \Sigma}$ the
{\em characteristic form} associated
to $u_i$, $i=1,2,3$.
Since
$
\chii_{j} (X_i) = \mu (X_i,X_j)= \delta_{ij},
$
each differential form $\om \in
\lau{\Om}{*}{}{ M\menos \Sigma}$ possesses a unique writing,
\begin{equation*}
\om = {{}_{_0}\omega} + \sum_{p=1}^3 \chii_p \wedge {{}_p\omega} + \sum_{1\leq p < q \leq 3}   \chii_{p} \wedge \chii_{q} \wedge {{}_{_{pq}}\omega}+
\chii_{1} \wedge \chii_{2} \wedge \chii_{3} \wedge {{}_{_{123}}\omega}  ,
\end{equation*}
where the coefficients ${{}_{_\bullet}\omega} $ are  {\em horizontal forms}, that is, they verify  $i_X \left({{}_{_\bullet}\omega} \right)=0$ for each $X \in \homo{\mathfrak{X}}{\Phi}{ M}$. This is the {\em canonical decomposition} of $\om$. For example,
$
d\beta  = {{}_{_0}(d\beta)} + \chii_1\wedge \ib{L}{X_1}\beta +
\chii_2\wedge  \ib{L}{X_2}\beta + \chii_3 \wedge  \ib{L}{X_3}\beta,
 $ for any  horizontal form $\beta $ (notice that this formula is no longer true if $\beta$ is not horizontal). Since $
\ib{L}{X_{i}}\chii_{j} = \chii_{[u_i,u_j]}$, with
$ 1\leq i,j\leq 3$, then
\begin{equation}
\label{lll}
\begin{array}{ll}
 \ib{L}{X_{1}}\chii_{1} = \ib{L}{X_{2}}\chii_{2}  = \ib{L}{X_{3}}\chii_{3}
    = 0  \hspace{2cm}&
   \ib{L}{X_{1}}\chii_{2} = - \ib{L}{X_{2}}\chii_{1}  =  \chii_{3}\\[,4cm]
 \ib{L}{X_{1}}\chii_{3} = - \ib{L}{X_{3}}\chii_{1}  =  -\chii_{2}&
    \ib{L}{X_{2}}\chii_{3} = - \ib{L}{X_{3}}\chii_{2}  = \chii_{1}
    \end{array}
    \end{equation}
    and we have the canonical decompositions
\be
  \label{dif}
      \left\{
    \begin{array}{rcl}\
d \chii_1 &= &e_1 - \chii_2 \wedge \chii_3 \\
d \chii_2 &= &e_2 + \chii_1 \wedge \chii_3 \\
d \chii_3 &=& e_3 - \chii_1 \wedge \chii_2.
\end{array}
\right.
\ee
Here, the forms $e_i$ are basic for $i=1,2,3$. Notice that $e_1-\chii_2\wedge\chii_3$ is the {\em Euler form} of the action of the maximal torus with fundamental vector field $X_1$, and that
$e_1^2+e_2^2+e_3^2$ is the  Euler form of the action $\Phi$ (see section \refp{section:morphisms}).

Consider $U \subset  M\menos \Sigma$ an equivariant open subset. If $\om \in \coho{\Om}{*}{ M\menos \Sigma,U}$ then the coefficients of its canonical decomposition  are horizontal forms of $\coho{\Om}{*}{ M\menos \Sigma,U}$.
The following Lemma is the key for the construction of the Gysin sequence.
Given an action  of $\ib{\Z}{2}$ on a vector space $E$  generated by the morphism $h \colon E \to E$, we shall write
\begin{equation}
\label{anti}
E^{-\ib{\Z}{2}} = \{ e \in E \tq h(e) = -e\},
\end{equation}
 the subspace of {\em antisymmetric elements}. Notice that $j \in \bi{\S}{3}$ acts naturally on $M^{\sbat}$.

\bl
\label{22}
\bee
\coho{H}{*}{\fracc{\lau{\BOm}{\cdot}{V}{M}}{\lau{\Om}{\cdot}{V}{M/\bi{\S}{3}} }} = \coho{H}{*-3}{M/\bi{\S}{3},\Sigma/\bi{\S}{3}}   \oplus  \left( \coho{H}{*-2}{M^{\sbat}} \right)^{-\ib{\Z}{2}}
\eee
\el
\pro
We consider the integration operator:
\begin{equation*}
\Int \colon\fracc{\lau{\BOm}{*}{V}{M}}{\lau{\Om}{*}{V}{M/\bi{\S}{3}} } \TO \lau{\Om}{*-3}{V}{M/\bi{\S}{3} , \Sigma/\bi{\S}{3}},
\end{equation*}
given by:
\begin{equation*}
\Int\left( < \om> \right) = (-1)^{\deg \om}  \ \ib{i}{X_3} \ib{i}{X_2} \ib{i}{X_1} \om.
\end{equation*}
 It is a well defined differential operator since

\begin{itemize}
\item[-] the tubular neighborhoods of the Thom-Mather's structure $\mathfrak{T} $ are invariant,
\item[-] the operator $\ib{i}{X_3} \ib{i}{X_2} \ib{i}{X_1}$ vanishes on $\Sigma$, and
\item[-] $\ib{i}{X}\ib{i}{X_3} \ib{i}{X_2} \ib{i}{X_1}\om =\ib{i}{X}d\ib{i}{X_3} \ib{i}{X_2} \ib{i}{X_1}\om =0$  for each $X \in \homo{\mathfrak{X}}{\Phi}{M}$\footnote{
$L_{A} i_{B} = i_{B} L_{A} +
 i_{[A,B]}, \quad  \forall A,B \in \mathfrak{X}(M).
$}.
\end{itemize}

Every form $\gamma \in \lau{\Om}{*-3}{V}{M/\bi{\S}{3} , \Sigma/\bi{\S}{3}}$ vanishes in  a neighborhood of $\Sigma$. So, the product $\chii_1 \wedge \chii_2 \wedge \chii_3 \wedge \gamma $ belongs to $\lau{\BOm}{*}{V}{M}$ (cf. \refp{lll}). Since $ \ib{i}{X_3} \ib{i}{X_2} \ib{i}{X_1} (\chii_1 \wedge \chii_2 \wedge \chii_3 \wedge  \gamma )=\gamma$ then we have the short exact sequence
\begin{equation}
\label{suc}
\xymatrix{
0\ar[r]  &
 \Ker^{*} \Int\
\ar@{^{(}->}[r] &\fracc{\lau{\BOm}{*}{V}{M}}{\lau{\Om}{*}{V}{M/\bi{\S}{3}} }  \ar[r]^-{\IInt} &
\lau{\Om}{*-3}{V}{M/\bi{\S}{3} , \Sigma/\bi{\S}{3}} \ar[r]& 0
}
\end{equation}

By Lemma     \ref{rela}, it suffices to prove the following:

\medskip

(a) $\coho{H}{*}{ \Ker^{*} \Int} =  \left( \coho{H}{*-2}{M^{\sbat}} \right)^{-\ib{\Z}{2}}$.

\medskip

(b) The associated connecting homomorphism $\delta $ vanishes.

\bigskip

\begin{center}
(a)
\end{center}

  For the sake of simplicity we put $\coho{\mathcal{A}}{*}{M} = \Ker^{*} \Int$. In fact we have
$
\coho{\mathcal{A}}{*}{M} = \fracc{\left\{ \om \in \lau{\BOm}{*}{V}{M} \tq \ib{i}{X_3} \ib{i}{X_2} \ib{i}{X_1}  \om =0\right\}}{\lau{\Om}{*}{V}{M/\bi{\S}{3}}}.
$
Analogously,  we define   $\coho{\mathcal{A}}{*}{M,\Sigma}$,  $\coho{\mathcal{A}}{*}{\Sigma}$ and  $\coho{\mathcal{A}}{*}{\Sigma,F} $. To get (a), it suffices to prove    the following facts:

\medskip

(a1) $\coho{H}{*}{\coho{\mathcal{A}}{*}{M}  } = \coho{H}{*}{\coho{\mathcal{A}}{*}{\Sigma} }$
 .

\medskip

(a2) $\coho{H}{*}{\coho{\mathcal{A}}{\cdot}{\Sigma} } = \coho{H}{*}{\coho{\mathcal{A}}{\cdot}{\Sigma,F}} $.

\medskip

(a3) $\coho{H}{*}{\coho{\mathcal{A}}{\cdot}{\Sigma,F} } =\left( \coho{H}{*-2}{M^{\sbat}} \right)^{-\ib{\Z}{2} } $.

\bigskip

\begin{center}
(a1)
\end{center}

Consider the inclusion
$L \colon  \coho{\mathcal{A}}{*}{M,\Sigma} \TO\coho{\mathcal{A}}{*}{M}
$
and the restriction
$
R \colon  \coho{\mathcal{A}}{*}{M} \to  \coho{\mathcal{A}}{*}{\Sigma},
$
which are differential morphisms. This gives the short sequence
\begin{equation*}
\xymatrix{
0\ar[r]  &
 \coho{\mathcal{A}}{*}{M,\Sigma}
\ar[r]^-L&\coho{\mathcal{A}}{*}{M}  \ar[r]^-R &
\coho{\mathcal{A}}{*}{\Sigma} \ar[r]& 0.
}
\end{equation*}
Notice that $R \rond L =0$. This short sequence is exact since:

\smallskip

$\bullet${\em The operator $R$ is an onto map.} Consider $\gamma \in \lau{\BOm}{*}{V}{\Sigma}$. We know that
  $\sigma (\gamma)\in \lau{\Om}{*}{V}{M}$ (cf. Lemma \ref{ext}). The result comes from:

$\rightsquigarrow$  $\sigma (\gamma)\in \lau{\BOm}{*}{V}{M}$.   Since $\tau_0$, $\tau_1$ are equivariant  and
$f  \rondp  \nu_0 $, $f  \rondp  \nu_1 $ are invariant.

\smallskip

$\rightsquigarrow$  $\ib{i}{X_3} \ib{i}{X_2} \ib{i}{X_1} \sigma (\gamma)=0$.  Since $\tau_0$, $\tau_1$ are equivariant and  $\rang \{X_1(x),X_2(x),X_3(X)\}\leq 2$  for any $x \in \Sigma $.

\smallskip

$\rightsquigarrow$  $R\left( <  \sigma (\gamma)  > \right) = < ( \sigma (\gamma))_1 > =  < \gamma  >$.

\smallskip

$\bullet$ {\em $\Ker R \subset \Ima L$.} Consider $\omega \in \lau{\BOm}{*}{V}{M}$ with $\ib{i}{X_3} \ib{i}{X_2} \ib{i}{X_1}  \om =0$ and $\ib{i}{X_j} \om_1 =0 \hbox{ for } j \in \{1,2,3\} $. Since $\tau_0$ and $\tau_1$ are equivariant  and $X_j=0$ on $F$ then $\ib{i}{X_j} \sigma (\om_1) =0 \hbox{ for }  j\in \{1,2,3\} $.
This gives $<  \sigma (\om_1) > =0$. Finally, we have $<  \om  > = <  \om -  \sigma (\om_1) > = L\left(  <  \om -  \sigma (\om_1) >\right)$
since $\left( \om - \sigma( \om_1) \right)_1 = \om_1 - ( \sigma( \om_1))_1 = \om_1 - \om_1 =0$.

\bigskip

Now, we will get (a1) by proving that  $\coho{H}{*}{\coho{\mathcal{A}}{\cdot}{M,\Sigma} } =0$.
By definition of Verona's  forms we have
$
\coho{\mathcal{A}}{*}{M,\Sigma} = \coho{\mathcal{A}}{*}{M,D}  \stackrel{excision}{=\!=\! = } \coho{\mathcal{A}}{*}{M\menos \Sigma,D \menos \Sigma},
$ where $D=D_0\cup D_1.$
 A straightforward calculation gives:
\begin{equation*}
\coho{H}{*}{\coho{\mathcal{A}}{\cdot}{M\menos \Sigma,D \menos \Sigma}} =
 \frac{
 \left\{ \om \in \lau{\BOm}{*}{}{M\menos \Sigma,D \menos \Sigma}\tq \ib{i}{X_3} \ib{i}{X_2} \ib{i}{X_1}  \om =0 \hbox{ and } \ib{i}{X_j} d\om =0 \hbox{ for  } j \in \{1,2,3\}\right\}
 }{
\lau{\Om}{*}{}{(M\menos \Sigma)/\bi{\S}{3},(D \menos \Sigma)/\bi{\S}{3}}
+
{\left\{ d\beta \tq \beta  \in \lau{\BOm}{*-1}{}{M\menos \Sigma,D \menos \Sigma} \hbox{ and } \ib{i}{X_3} \ib{i}{X_2} \ib{i}{X_1}  \beta =0 \right\}
}}
\end{equation*}
Let $\om $ be a differential form of $\lau{\BOm}{*}{}{M\menos \Sigma,D \menos \Sigma}$ verifying $\ib{i}{X_3} \ib{i}{X_2} \ib{i}{X_1}  \om =0$ and $\ib{i}{X_j} d\om =0 \hbox{ for  } j\in \{1,2,3\}$. Then
\begin{equation*}
\begin{array}{ccc}
&&d\underbrace{\left( \chii_1 \wedge i_{X_3} i_{X_2} \om - \chii_2 \wedge i_{X_3} i_{X_1} \om
+ \chii_3 \wedge i_{X_2} i_{X_1} \om\right)}_{\beta} \\ \om\hspace{1cm}&=\hspace{2cm}&+\\&&
\underbrace{- e_1 \wedge i_{X_3} i_{X_2} \om + e_2 \wedge i_{X_3} i_{X_1} \om
- e_3 \wedge i_{X_2} i_{X_1} \om  + {{}_{_0}\om}}_{\alpha}
 \end{array}
 \end{equation*}
(cf. \refp{dif}) with $\beta \in  \lau{\BOm}{*-1}{}{M\menos \Sigma,D}$, verifying $\ib{i}{X_3} \ib{i}{X_2} \ib{i}{X_1}\beta =0 $, and $\alpha \in \lau{\Om}{*}{}{(M\menos \Sigma)/\bi{\S}{3},D/\bi{\S}{3}} $. This implies $\coho{H}{*}{\coho{\mathcal{A}}{\cdot}{M\menos \Sigma,D \menos \Sigma} } =0$ and then $\coho{H}{*}{\coho{\mathcal{A}}{\cdot}{M,\Sigma} } =0$.

\bigskip

\begin{center}
(a2)
\end{center}

Consider the inclusion
$L \colon  \coho{\mathcal{A}}{*}{\Sigma,F} \hookrightarrow \coho{\mathcal{A}}{*}{\Sigma}$
which is a differential morphism.
It suffices to prove that $L$ is an onto map.

Let us consider a smooth function $f \colon ]0,\infty[ \to [0,1]$ verifying $f \equiv 0$ on $[3,\infty[$ and $f \equiv 1$ on $]0,2]$. Notice that the composition $ f  \rondp  \nu_0 \colon M \to [0,1]$ is a  smooth invariant map.
So, for each $\gamma \in \coho{\Om}{*}{F}$ we have $ \sigma (\gamma) =( f  \rondp  \nu_0)\tau^*_0  \gamma \in \lau{\Om}{*}{}{M}$. This differential form verifies

\smallskip

$\rightsquigarrow$ $\sigma (\gamma) \in \lau{\Om}{*}{V}{\Sigma}$. Since  $(f \rondp \nu_0) \equiv 1 $ on $\Delta_0$ then $\sigma (\gamma) =\tau_0^* \gamma$ on  $\Delta_0 \menos F$. This gives
$(\sigma_0 (\gamma))_0 = \gamma$.

\smallskip

$\rightsquigarrow$  $\sigma (\gamma)\in \lau{\BOm}{*}{V}{\Sigma}$.   Since $\tau_0$ is an equivariant map and $f  \rondp  \nu_0 $ is an invariant map.

 \smallskip

$\rightsquigarrow$ $\ib{i}{X_j} \sigma (\gamma) =0 \hbox{ for }  j\in \{1,2,3\} $ since $\tau_0$ is an equivariant map and $X_j=0$ on $F$.

\nt Then
 $<\sigma (\gamma) > =0$ on $\coho{\mathcal{A}}{*}{\Sigma}$.

Let $<\om>$ be a class of $\coho{\mathcal{A}}{*}{\Sigma}$. We can write:
$
<\om> = <\om- \sigma ((\omega)_0)> =  L\left(  <  \om -   \sigma ((\omega)_0)>\right)$
since $\left( \om -  \sigma ((\omega)_0)\right)_0= \om_0 - ( \sigma( \om_0))_0 = \om_0 - \om_0 =0$.  This proves that $L$ is an onto map.

\bigskip

\begin{center}
(a3)
\end{center}

By definition of Verona's differential forms we have
\begin{equation*}
\coho{\mathcal{A}}{*}{\Sigma,F} = \coho{\mathcal{A}}{*}{\Sigma,\Delta_0}  \stackrel{excision}{=\!=\! = } \coho{\mathcal{A}}{*}{\Sigma\menos F,\Delta_0 \menos F}
=
\frac{ \lau{\BOm}{*}{}{\Sigma\menos F,\Delta_0 \menos F}}{\lau{\Om}{*}{}{(\Sigma\menos F)/\bi{\S}{3},(\Delta_0 \menos F)/\bi{\S}{3}}} .\end{equation*}
The isotropy subgroup of a point of $\Sigma\menos F$ is conjugated to $\sbat$ or $N(\sbat)$ (cf. \cite[Th. 8.5, pag. 153]{Br}). We consider the manifold
$
\Gamma = \left(\Sigma\menos F\right)^{\sbat}.
$
A straightforward calculation gives that $\Sigma\menos F$ is $G$-equivariant diffeomorphic to
\bee
 \bi{\S}{3}\times_{N(\sbat)} \Gamma= \left( \bi{\S}{3}/\sbat\right) \times_{N(\sbat)/\sbat} \Gamma=
 \bi{\S}{2} \times_{\ib{\Z}{2}}   \Gamma.
\eee
Notice that  $  \Gamma/\ib{\Z}{2} =(\Sigma\menos F) /\bi{\S}{3}$.  Let $\Gamma_0 $ be the open subset $\Gamma \cap \Delta_0$ of $\Gamma$. Analogously we have $\Delta_0 \menos F = \bi{\S}{2}\times_{\ib{\Z}{2}} \Gamma_0$ and  $\Gamma_0/\ib{\Z}{2} = \left(\Delta_0 \menos F\right)/\bi{\S}{3}$.

The $\ib{\Z}{2}$-action on $ \bi{\S}{2}$ is generated by
$(x_0,x_1,x_2) \mapsto (-x_0,-x_1,-x_2)$\footnote{This map is
induced by $j \colon \bi{\S}{3} \to \bi{\S}{3}$ defined by $j(u) = u
\cdot j$ (see \cite[Example 17.23]{BT}).}. Then, the
$\ib{\Z}{2}$-action on $\coho{H}{0}{\bi{\S}{2}}$ (resp.
$\coho{H}{2}{\bi{\S}{2}}$) is the identity $\Ide$ (resp. $-\Ide$).
The $\ib{\Z}{2}$-action on $  \Gamma$  is induced by $\Phi(j,-)$.
The K\"{u}nneth formula gives
\begin{eqnarray*}
 \coho{H}{*}{
  \lau{\BOm}{*}{}{\Sigma\menos F,\Delta_0 \menos F}}
&=&
\coho{H}{*}{
  \lau{\BOm}{\cdot}{}{ \bi{\S}{2} \times_{\ib{\Z}{2}}   \Gamma,\bi{\S}{2} \times_{\ib{\Z}{2}}   \Gamma_0}
 }
 =
 \coho{H}{*}{
  \lau{\BOm}{\cdot}{}{ \bi{\S}{2} \times  \Gamma,\bi{\S}{2} \times   \Gamma_0}^{\ib{\Z}{2}}
 }\\
 & =&
 \coho{H}{*}{
  \lau{\BOm}{\cdot}{}{ \bi{\S}{2} \times  \Gamma,\bi{\S}{2} \times   \Gamma_0}
 }^{\ib{\Z}{2}}
=
\left(
\lau{H}{*}{}{\bi{\S}{2}} \otimes   \lau{H}{*}{}{  \Gamma,   \Gamma_0
 }\right)^{\ib{\Z}{2}} \\
 &= &
\left(\lau{H}{0}{}{\bi{\S}{2}} \otimes  \lau{H}{*}{}{  \Gamma,   \Gamma_0
 }\right)^{\ib{\Z}{2}}
 \oplus
\left(\lau{H}{2}{}{\bi{\S}{2}} \otimes   \lau{H}{*-2}{}{  \Gamma,   \Gamma_0
 }\right)^{\ib{\Z}{2}}
=
 \left( \lau{H}{*}{}{  \Gamma,   \Gamma_0
 }\right)^{\ib{\Z}{2}}
 \oplus
 \left( \lau{H}{*-2}{}{  \Gamma,   \Gamma_0
 }\right)^{-\ib{\Z}{2}}
  \\
&= &
 \lau{H}{*}{}{  \Gamma/\ib{\Z}{2},   \Gamma_0/\ib{\Z}{2} }
 \oplus
 \left( \lau{H}{*-2}{}{  \Gamma,   \Gamma_0
 }\right)^{-\ib{\Z}{2}} =
 \lau{H}{*}{}{(\Sigma \menos F)/\bi{\S}{3}, (\Delta_0\menos F)/\bi{\S}{3} }
 \oplus
 \left( \lau{H}{*-2}{}{  \Gamma,   \Gamma_0
 }\right)^{-\ib{\Z}{2}} ,
 \end{eqnarray*}
 and then
\begin{eqnarray*}
\coho{H}{*}{\coho{\mathcal{A}}{\cdot}{\Sigma\menos F,\Delta_0 \menos F }}  &=&
 \coho{H}{*}{\frac{ \lau{\BOm}{\cdot}{}{\Sigma\menos F,\Delta_0 \menos F}}{\lau{\Om}{\cdot}{}{(\Sigma\menos F)/\bi{\S}{3},(\Delta_0 \menos F)/\bi{\S}{3}}}}
=
\left( \lau{H}{*-2}{}{  \Gamma,   \Gamma_0
 }\right)^{-\ib{\Z}{2}} =
 \left(\lau{H}{*-2}{}{(\Sigma\menos F) ^{\sbat}, (\Delta_0 \menos F) ^{\sbat}}
 \right)^{-\ib{\Z}{2}}\\
 & \stackrel{excision}{=\!=\! = } &
 \left(\lau{H}{*-2}{}{\Sigma^{\sbat}, \Delta_0^{\sbat}}
 \right)^{-\ib{\Z}{2}} \stackrel{retraction}{=\!=\! = }
 \left(\lau{H}{*-2}{}{\Sigma^{\sbat}, F^{\sbat}}
 \right)^{-\ib{\Z}{2}}.
 \end{eqnarray*}
Consider  the long exact sequence associated to the $\ib{\Z}{2}$-invariant pair $\left(\Sigma^{\sbat}, F^{\sbat}\right)$:
\begin{equation*}
\cdots \to
\left( \lau{H}{i-1}{}{F ^{\sbat}}\right)^{- \ib{\Z}{2} } \to
\left( \lau{H}{i}{}{\Sigma^{\sbat}, F^{\sbat} }\right)^{- \ib{\Z}{2}}
 \to
\left( \lau{H}{i}{}{\Sigma^{\sbat}}\right)^{- \ib{\Z}{2}}
 \to
\left( \lau{H}{i}{}{F^{\sbat}}\right)^{- \ib{\Z}{2}}
 \to
\cdots .
\end{equation*}
Since the action of $ \ib{\Z}{2}$ on $F^{\sbat}=F$ is  trivial, then $\left( \lau{H}{i}{}{F ^{\sbat}}\right)^{- \ib{\Z}{2}} =0$. On the other hand, we have $\Sigma^{\sbat} = M^{\sbat}$. This gives
$
  \left(\lau{H}{*-2}{}{\Sigma^{\sbat}, F ^{\sbat}}
 \right)^{-\ib{\Z}{2}} =   \left(\lau{H}{*-2}{}{M^{\sbat}}
 \right)^{-\ib{\Z}{2}}. $

 \bigskip

\begin{center}
(b)
\end{center}
Notice that the connecting morphism $\delta$ is defined by $\delta ([\zeta]) =\pm \left[ <  d (\chii_1 \wedge \chii_2 \wedge \chii_3) \wedge \zeta > \right]$. We have $\delta \equiv 0$ since $\zeta_1 =0$ (cf (a1)).
\qed

\bt
\label{gysin}
Given any smooth action $\Phi \colon \bi{\S}{3} \times M \TO M$ we have the Gysin sequence
\begin{equation*}
\xymatrix{
\cdots \ar[r] &
\coho{H}{i}{M}\ar[r] &
\coho{H}{i-3}{M/\bi{\S}{3},\Sigma/\bi{\S}{3}}   \oplus  \left( \coho{H}{i-2}{M^{\sbat}} \right)^{-\ib{\Z}{2}}
\ar[r] &\coho{H}{i+1}{M/\bi{\S}{3}} \ar[r]
 &  \coho{H}{i+1}{M}\ar[r] & \cdots &
}
\end{equation*}
where
$\Sigma$ is the subset of  points of $M$ whose isotropy group is infinite, the $\ib{\Z}{2}$-action is induced by $j \in \bi{\S}{3}$ and $(-)^{-\ib{\Z}{2}}$ denotes the subspace of antisymmetric elements.
\et
\pro
Consider the short exact sequence
\begin{equation}
\label{cortogysin}
\xymatrix
{
0 \ar[r]&
\lau{\Om}{*}{V}{M/\bi{\S}{3}} \ar@{->}[r]&
\lau{\BOm}{*}{V}{M}\ar[r] &
\fracc{\lau{\BOm}{*}{V}{M}}{\lau{\Om}{*}{V}{M/\bi{\S}{3}} }
\ar[r] &0,
}
\end{equation}
take its associated long exact sequence and then, apply Lemma \ref{rela},  \refp{ver} and Lemma \ref{22}.
\qed

\prgg{Example}
Consider the connected sum $M=\C\P^2\ \#\
\C\P^2\cong\left(\S^3\times[0,1]\right)/\sim$, with
$$
((z_1,z_2),i)\sim ((z\cdot z_1,z\cdot z_2),i),\quad i=0,1,
$$
for all $z\in\S^1$ and $(z_1,z_2)\in\S^3$ in complex coordinates. The product of $\S^3$ induces
on $M$ the action:
$$
g\cdot[h,t]=[g\cdot h,t],\quad \forall g,h\in\S^3, \forall
t\in[0,1].
$$
For this action, we have:
$$
\Sigma =\left(\S^3\times\{0,1\}\right)/\sim\ \cong\S^2\times
\{0,1\},\qquad F=\emptyset,
$$
$$
M/\S^3\cong[0,1],\quad \Sigma/\S^3\cong\{0,1\},\quad M^{\S^1}\cong
\{N,S\}\times\{0,1\},
$$
where $N$ and $S$ stand for the North and South poles of $\S^2$. The $\Z_2$- action on $M^{\S^1}$
is determined by $j\in\S^3$, which induces the antipodal map on $\S^2$, and so, interchanges
its poles. Thus, the exotic term that appears in the central part of the Gysin Sequence is not trivial:
$$
\lau{H}{2}{}{M} \stackrel{\cong}{\TO}
\left( \coho{H}{0}{M^{\sbat}} \right)^{-\ib{\Z}{2}}
 =
 \left( \coho{H}{0}{\{N,S\}\times\{0,1\}} \right)^{-\ib{\Z}{2}}
 \cong
 \R\oplus\R.
$$

\prgg{Morphisms}\label{section:morphisms} We describe the morphisms of the Gysin sequence.

\begin{center}
$\textcircled{\tiny  1}\colon \coho{H}{*}{M/\bi{\S}{3}} \TO  \coho{H}{*}{M} $
\end{center}
It is the pull-back $\pi^*$ of the canonical projection $\pi \colon M \to M/\bi{\S}{3}$ (cf.  Lemma \ref{rela}).

\begin{center}
$\textcircled{\tiny  2}\colon \coho{H}{*}{M}\TO
\coho{H}{*-3}{M/\bi{\S}{3},\Sigma/\bi{\S}{3}}   \oplus  \left( \coho{H}{*-2}{M^{\sbat} } \right)^{-\ib{\Z}{2}}$
\end{center}

We have already seen that the first component of this morphism is induced by $\Int_{\bi{\S}{3}} [\om] = [\ib{i}{X_3} \ib{i}{X_2} \ib{i}{X_1}\om]$.
 For the second component we keep track of the isomorphisms given by Lemma \ref{22} and we get that it is defined by:
$
[\om] \mapsto \class \left(\Int_{\bi{\S}{2}} \left( \om_1 - \sigma (\iota^*\om_1)\right) \right).$
\begin{center}
$\textcircled{\tiny  3}\colon
\coho{H}{*-3}{M/\bi{\S}{3},\Sigma/\bi{\S}{3}}   \oplus  \left( \coho{H}{*-2}{M^{\sbat} } \right)^{-\ib{\Z}{2}} \TO \coho{H}{*+1}{M/\bi{\S}{3}}$
\end{center}
A straightforward calculation using sequences \refp{suc} and \refp{cortogysin} gives that the connecting morphism $\textcircled{\tiny  3}$ of the Gysin sequence  sends:
\begin{itemize}
\item $[\zeta] \in \coho{H}{*-3}{M/\bi{\S}{3},\Sigma/\bi{\S}{3}} $ to -$\left[\left( e_1^2 + e_2^2 + e_3^2 \right) \wedge \zeta\right]$, and
\item $[\xi] \in  \left( \coho{H}{*-2}{M^{\sbat} } \right)^{-\ib{\Z}{2}} =  \left(\lau{H}{*-2}{}{\Sigma^{\sbat}, F ^{\sbat}}
 \right)^{-\ib{\Z}{2}}$ to
$
[d\sigma \wedge \epsilon \wedge \tau_1^* \xi]
$
where $\epsilon $ is an Euler form of the restriction $\Phi_1 \colon \sbat \times \left( \tau_1^{-1}\left( \Sigma^{\sbat}\right) \menos \Sigma^{\sbat}\right) \to\left( \tau_1^{-1}\left( \Sigma^{\sbat}\right) \menos \Sigma^{\sbat}\right)$ of $\Phi$.
\end{itemize}
Since $e_1^2 + e_2^2 + e_3^2$ is not a Verona's form, then it does not define a class of $\coho{H}{4}{M/\bi{\S}{3}}$. Nevertheless, it does generate a class in  the intersection cohomology group $\lau{\IH}{4}{\per{4}}{M/\bi{\S}{3}}$ (as in the semi-free case of \cite{S3}).

\prgg{Remarks}

\Zati We have {
$ \left( \coho{H}{*}{M^{\sbat} } \right)^{-\ib{\Z}{2}} = \coho{H}{*}{M^{\sbat}}\Big/  \coho{H}{*}{M^{\sbat}/\ib{\Z}{2}}
 $. Let us see that. The correspondence $\om \mapsto \left( \fracc{\om + j^* \om}{2}, \fracc{\om - j^* \om}{2}\right)$ establishes the isomorphism
$ \coho{\Om}{*}{M^{\sbat} }  =
 \left( \coho{\Om}{*}{M^{\sbat} } \right)^{\ib{\Z}{2}} \oplus  \left( \coho{\Om}{*}{M^{\sbat} } \right)^{- \ib{\Z}{2}}
 =
\coho{\Om}{*}{M^{\sbat}/\ib{\Z}{2}}  \oplus  \left( \coho{\Om}{*}{M^{\sbat} } \right)^{- \ib{\Z}{2}}
$ and hence,
$
 \coho{H}{*}{M^{\sbat} }  =
\coho{H}{*}{M^{\sbat}/\ib{\Z}{2}}  \oplus  \left( \coho{H}{*}{M^{\sbat} } \right)^{- \ib{\Z}{2}}.
$
This gives the claim.

\zati  Let us suppose that the action is semi-free, almost free or free.
Then,  $j$ acts trivially on $M^{\sbat} =F$, and hence,
 we have a long exact sequence
\begin{equation*}
\cdots \to \coho{H}{i}{M} \to \coho{H}{i-3}{M/\bi{\S}{3},F}\to \coho{H}{i+1}{M/\bi{\S}{3}} \to \coho{H}{i+1}{M} \to \cdots.
\end{equation*}

\zati Let us suppose that there is not a point of $M$ whose isotropy subgroup is conjugated to $\sbat$.
Then, we have a long exact sequence
\begin{equation*}
\cdots \to \coho{H}{i}{M} \to \coho{H}{i-3}{M/\bi{\S}{3},\Sigma/\bi{\S}{3}}\to \coho{H}{i+1}{M/\bi{\S}{3}} \to \coho{H}{i+1}{M} \to \cdots .
\end{equation*}
since $j$ acts trivially on $M^{\sbat}=\left\{ x \in M \tq \bi{\S}{3}_x = \S^3 \hbox{ or } N(\sbat)\right\}$.

\prgg{Actions over $\sbat$}

Using the Gysin sequence we have constructed, we now give a list of
all the different cohomologies of a $\bi{\S}{3}$-manifold $M$ having
the circle as orbit space\footnote{In fact, we give  the Poincar\'e
polynomial $\ib{P}{M}$  of $M$.}. By geometrical reasons, the orbit
space is composed by just one stratum, the whole circle. Following
the nature of the orbits, we distinguish four  cases.

\Zati All orbits are of dimension 3. We have $\ib{P}{M} =
1+t+t^3+t^4$. This is the case of the manifold $\bi{\S}{3} \times
\sbat$, where $\bi{\S}{3}$ acts by multiplication on the left
factor.

\zati All orbits are isomorphic to $\bi{\S}{2}$. We distinguish two
cases following wether the covering $M^{\sbat} \to M^{\sbat}
/\ib{\Z}{2} = M/M^{\sbat} $ is trivial or not. In the first case we
have $\ib{P}{M} = 1+t+t^2+t^3$. This is the case of the manifold
$\bi{\S}{2} \times \sbat$, where $\bi{\S}{3}$ acts by multiplication
on the left factor. In the second case we have $\ib{P}{M} = 1+t$, as
is the case of the manifold $\bi{\S}{2} \ib{\times}{\bi{\Z}{2} }
\sbat$ where $\bi{\S}{3}$ acts by multiplication on the left factor.

\zati All orbits are isomorphic to $\bi{\R\P}{2}$. In this case, we
have $\ib{P}{M} = 1+t$. This is the case of the manifold
$\bi{\R\P}{2}\times \sbat$ where $\bi{\S}{3}$ acts by multiplication
on the left  factor.

\zati All orbits are points. We have $\ib{P}{M} = 1+t$. This
corresponds to the manifold $\sbat$ where $\bi{\S}{3}$ acts
ineffectively.

\end{document}